\documentclass[10pt,twoside]{article}
\usepackage{amsmath}
\usepackage{amsfonts}

\newcommand{\sm}{\setminus}
\newcommand{\szego}{Szeg\"o }

\newcommand{\Si}{\Sigma}
\newcommand{\inv}{^{-1}}
\newcommand{\kahler}{K\"ahler }
\newcommand{\sqrtn}{\sqrt{N}}
\newcommand{\wt}{\widetilde}

\newcommand{\PP}{{\mathbb P}}
\newcommand{\N}{{\mathbb N}}
\newcommand{\R}{{\mathbb R}}
\newcommand{\C}{{\mathbb C}}

\newcommand{\Z}{{\mathbb Z}}

\newcommand{\CP}{\C\PP}

\newcommand{\E}{{\mathbf E}}

\newcommand{\half}{{\frac{1}{2}}}
\newcommand{\vol}{{\operatorname{Vol}}}

\newcommand{\SU}{{\operatorname{SU}}}
\newcommand{\FS}{{{\operatorname{FS}}}}

\renewcommand{\phi}{\varphi}

\newcommand{\ocal}{\mathcal{O}}

\newcommand{\ga}{\gamma}

\newcommand{\la}{\lambda}

\newcommand{\de}{\delta}

\newcommand{\om}{\omega}

\newtheorem{theo}{Theorem}

\usepackage{Latex-document}

\title{\bf  Asymptotics of Polynomials and Eigenfunctions \vskip -2mm\vskip 6mm}
\author{S. Zelditch\vspace*{-0.5cm}\thanks{Department of Mathematics,
Johns Hopkins University, Baltimore, Maryland 21218, USA. E-mail:
zelditch@math.jhu.edu}}
\date{\vspace{-8mm}}

\begin{document}

\maketitle

\thispagestyle{first} \setcounter{page}{733}

\begin{abstract}

\vskip 3mm

We review some recent results on asymptotic properties of
polynomials  of large degree, of general holomorphic sections of
high powers of  positive line bundles over \kahler manifolds, and
of Laplace eigenfunctions of large eigenvalue on compact
Riemannian manifolds. We describe statistical patterns in the
zeros, critical points and $L^p$ norms of random polynomials and
holomorphic sections, and the influence of
 the Newton polytope on these patterns. For eigenfunctions, we discuss $L^p$ norms and mass concentration of individual eigenfunctions and their
relation to dynamics of the geodesic flow.

\vskip 4.5mm


\end{abstract}

\vskip 12mm

\section{Introduction}

\vskip-5mm \hspace{5mm}

 In many measures  of  `complexity',    eigenfunctions
$ \sqrt{\Delta} \phi_{\lambda} = \lambda \phi_{\lambda}$ of first
order elliptic  operators  behave like polynomials $p(x) =
\sum_{|\alpha| \leq x} c_{\alpha} x^{\alpha} $
 of  degree $N \sim \lambda$ \cite{DF}. The basic example we have
in mind is the Laplacian  $\Delta$  on a compact Riemannian
manifold $(M, g)$, but the same is true of Schroedinger operators.
The comparison is more than an analogy, since polynomials of
degree $N$ are eigenfunctions of a first order elliptic system.

 The comparison between  eigenfunctions and polynomials
 is an essentially  local one,  most accurate
 on  small balls
$B(x_0, \frac{1}{\lambda})$. Globally,   eigenfunctions reflect
the dynamics of the geodesic flow $G^t$ on the unit (co-)tangent
bundle $S^*M$. This is one of the principal themes  of quantum
chaos.

In this article, we review some recent results on the asymptotics
of polynomials and eigenfunctions, concentrating on our work  in
collaboration with P. Bleher, A. Hassell, B. Shiffman, C. Sogge,
J. Toth and M. Zworski. A unifying feature is the asymptotic
properties of reproducing kernels, namely \szego kernels
$\Pi_N(z,w)$ in the case of polynomials, and spectral projections
$E_{\lambda}(x, y)$ for intervals $[\lambda, \lambda + 1]$ in the
case of eigenfunctions of $\sqrt{\Delta}$.  For other recent
expository articles, see \cite{JNT, Z}.

\section{Polynomials}

\vskip-5mm \hspace{5mm }

There are several sources of interest in random polynomials. One
is  the desire to understand typical properties of real and
complex algebraic varieties, and how they depend on the
coefficients of the defining equations. Another is their use as a
model for the local behavior of more general eigenfunctions. A
third is that they may be viewed  as the eigenvectors of  random
matrices. Just as random matrices model the spectra of `quantum
chaotic' systems, so random polynomials model their
eigenfunctions.

\subsection{{\boldmath $SU(m + 1)$} polynomials on {\boldmath $\CP^m$} and holomorphic sections}
\vskip-5mm \hspace{5mm }

Complex polynomials of degree $\leq p$ in $m$  variables form the
vector  space
$${\mathcal P}_p^m: = \{f(z_1, \dots, z_m) =
\sum_{\alpha \in \N^m: |\alpha| \leq p} c_{\alpha} z_1^{\alpha_1}
\cdots z_m^{\alpha_m},\;\ c_{\alpha} \in \C\}\;.$$ To put a
probability measure on ${\mathcal P}_p^m$ is to regard the
coefficients $c_{\alpha}$ as random variables. The simplest
measures are Gaussian measures corresponding to inner products on
${\mathcal P}_p^m$. By homogenizing $f$ to $ F(z_0, z_1, \cdots,
z_m)$ of degree $p$, we may identify  ${\mathcal P}_p^m$  with the
space
   $H^0(\CP^m, {\mathcal
O}(p))$ of holomorphic sections of the $p$th power of the
hyperplane bundle. It carries  the standard $\SU(m+1)$-invariant
Fubini-Study inner product $\langle F_1, \bar F_2\rangle_{FS} =
\int_{S^{2m+1}} F_1 \bar F_2\, d\sigma\;,$ where $d\sigma$ is Haar
measure on the $(2m+1)$-sphere $S^{2m+1}$. An orthonormal basis of
$H^0(\CP^m, {\mathcal O}(p))$ is given by
$\{\frac{z^{\alpha}}{||z^{\alpha}||_{FS}}\}.$ The corresponding
$\SU(m+1)$-invariant  Gaussian measure $\gamma_\de$ is defined by
$$ d
\gamma_\de (s) = \frac{1}{\pi^{k_p}}e^{-|\lambda|^2} d\lambda,\;\;
s = \sum_{|\alpha| \leq p} \lambda_{\alpha}
\frac{z^{\alpha}}{||z^{\alpha}||_{FS}}.$$ Thus, the coefficients
$\lambda_{\alpha}$ are independent complex Gaussian random
variables with mean zero and variance one.

More generally, we can define Gaussian ensembles of holomorphic
sections $H^0(M, L^N)$ of powers of a positive line bundle over
any \kahler manifold $(M, \omega)$. Endowing   $L$ with the unique
hermitian metric $h$ of curvature form $\omega$, we induce an
inner product $\langle, \rangle$ on $H^0(M, L^N)$ and a Gaussian
measure $\gamma_N$. We denote the unit sphere in $H^0(M, L^N)$
relative to this inner product by $SH^0(M, L^N)$. The Haar measure
on $SH^0(M, L^N)$ will be denoted $\mu_N$. It is closely related
to the Gaussian measure.

\subsection{Zeros}
\vskip-5mm \hspace{5mm }

The problems we  discuss in this section involve the geometry of
zeros of sections $s \in H^0(M, L^N)$ of general positive line
bundles. There is a similar story for critical points.

\begin{itemize}

\item {\bf Problem 1} How are the simultaneous zeros  $Z_s = \{z: s_1(z) = \cdots = s_k(z) = 0\}$ of a k-tuple $s=(s_1,\dots,s_k)$ of typical  holomorphic sections
distributed?

\item {\bf Problem 2} How are the zeros correlated? When $k = m$, the simultaneous
zeros form a discrete set. Do  zeros repel each other like charged
particles?  Or behave independently
 like particles of an ideal gas? Or attract like gravitating particles?

\end{itemize}

By the distribution of zeros we mean either the current of
integration over  $Z_s$ or more simply   the Riemannian
$(2m-2k)$-volume measure $(|Z_s|,\phi)=\int_{Z_s}\phi
d\vol_{2m-2k}\,.$ By  the $n$-point zero correlation
 functions, we mean
the  generalized functions
$$K_{nk}^N(z^1,\dots,z^n)dz\, = \E |Z_s|^n,$$
where $|Z_s|^n$ denotes  the product of the  measures $|Z_s|$ on
 the punctured product
$M_n=\{(z^1,\dots,z^n) \in M \times \cdots \times M: z^p\ne z^q \
\ {\rm for} \ p\ne q\}$ and where $dz$ denotes the product volume
form on $M_n$.

The answer to Problem 1  is that zeros almost surely become
uniformly distributed relative to the curvature $\omega$ of the
line bundle \cite{SZ1}. Curvature causes sections to oscillate
more rapidly and hence to vanish more often. More precisely, we
consider the space ${\mathcal S} = \prod_{N=1}^{\infty} SH^0(M,
L^{N})$ of random sequences, equipped with the product measure
measure $\mu = \prod_{N=1}^{\infty} \mu_N$.
 An element in ${\mathcal S}$ will be denoted ${\bf s} = \{s_N\}$. Then,
$\frac{1}{N} Z_s \to \omega, $ as $N \to \infty$ for almost every
${\bf s}.$

The answer to Problem 2 is more subtle: it depends on the
dimension. We assume $k =m$ so that almost surely the simultaneous
zeros of the $k$-tuple of sections form a discrete set. We find
that these zeros  behave almost independently if they are of
distance $\geq \frac{D}{\sqrt{N}}$ apart for $D \gg 1$.  So they
only interact on distance scales of size $\frac{1}{\sqrt{N}}$.
Since also the  density of zeros in a unit ball $B_{1}(z_0)$
around $z_0$  grows like $ N^m$, we rescale the zeros in the
$1/\sqrt{N}$-ball $B_{1/\sqrt{N}}(z_0)$ by a factor of $\sqrt{N}$
to get configurations of zeros with a constant density as $N \to
\infty.$ We thus rescale the correlation functions and take the
{\it scaling limits}
  \begin{equation}\label{slcd}\wt K^\infty_{nkm}(z^1,\dots,z^n)= \lim_{N\to\infty}K^N_{1k}(z_0)^{-n}
K_{nk}^N(z_0+\frac{z^1}{\sqrtn},\dots,z_0+\frac{z^n}{\sqrtn})\,.\end{equation}
In \cite{BSZ1}, we proved that the scaling limits of these
correlation functions were universal, i.e. independent of $M, L,
\omega, h$. They depend
 only on the dimension
$m$ of the manifold and the codimension $k$ of the zero set.

In \cite{BSZ2}, we found explicit formulae for these universal
scaling limits.
 In  the case $n = 2$,   $\wt K_{2km}^\infty(z^1,z^2)$, depends only on the
distance between the points $z^1,z^2$, since it is universal and
hence invariant under rigid motions.  Hence it may be written as:
\begin{equation} \wt K_{2km}^\infty(z^1,z^2) = \kappa_{km}(|z^1-z^2|)\,.
\end{equation}
We refer to \cite{BSZ1} for details.

\begin{theo} \label{2mm} {\rm \cite{BSZ2}}  The pair correlation functions of zeros when $k = m$ are given by
\begin{equation} \label{leading} \kappa_{mm}(r)= \left\{ \begin{array}{ll}\frac{m+1}{4}
r^{4-2m} + O(r^{8-2m})\,, & \qquad\mbox{as }\ r\to 0\, \\ [2mm]
1 + O(e^{- C r^2}), \;\; (C > 0) & \qquad\mbox{as }
r \to \infty.
\end{array} \right.\end{equation}
\end{theo}

 When
$m = 1, \kappa_{mm}(r) \to 0$ as $r \to 0$ and one has ``zero
repulsion.'' When $m = 2$, $\kappa_{mm}(r) \to  3/4$ as $r \to 0$
and one has a kind of neutrality. With $m \geq 3$, $\kappa_{mm}(r)
\nearrow \infty$ as $r \to 0$ and there is some kind of attraction
between zeros. More precisely, in dimensions greater than 2, one
is more likely to find a zero at a small distance $r$ from another
zero than at a small distance $r$ from a given point; i.e., zeros
tend to clump together in high dimensions.

 One can understand this dimensional dependence
heuristically in terms of the geometry   of the discriminant
varieties ${\mathcal D}_N^m \subset H^0(M, L^N)^m$ of systems $S =
(s_1, \dots, s_m)$ of $m$ sections with a `double zero'. The
`separation number' $sep(F)$ of a system is the minimal distance
between a pair of its zeros. Since the nearest element of
${\mathcal D}_N^m $ to $F$ is likely to have a simple double zero,
one expects: $sep(F) \sim \sqrt{dist(F, {\mathcal D}_N^m })$.
Now,the degree of ${\mathcal D}_N^m$ is approximately $N^m$.
Hence, the tube $({\mathcal D}_N^m )_{\epsilon}$ of radius
$\epsilon$ contains a volume $\sim \epsilon^2 N^m$. When $\epsilon
\sim N^{-m/2}$, the tube  should cover $P H^0(M, L^N)$. Hence, any
section should have a pair of zeros whose separation is  $\sim
N^{-m/4}$ apart. It is clear that this separation is larger than,
equal to or less than $N^{-1/2}$ accordingly as $m = 1, m = 2, m
\geq 3.$

\subsection{Bergman-\szego kernels}
\vskip-5mm \hspace{5mm }

A key object in the proof of these results is the Bergman-\szego
kernel $\Pi_{N}(x, y)$, i.e. the kernel of the orthogonal
projection onto $H^0(M, L^N)$ with respect to the \kahler form
$\omega.$ For instance, the expected distribution of zeros is
given by $\E_{N} (Z_f) = \frac{\sqrt{-1}}{2\pi}
\partial  \bar{\partial} \log \Pi_{N}(z,z) + \omega$.
Of even greater use is  the  joint probability distribution (JPD)
$D_N(x^1,\dots,x^n;\xi^1,\dots,\xi^n;z^1,\dots,z^n)$
 of the random variables
$x^j(s) =s(z^j), \ \xi^j(s) =\nabla s(z^j),$ which may be
expressed in terms of $\Pi_N$ and its derivatives. In turn, the
correlation functions may be expressed in terms of the JPD by
 $K^N(z^1,\dots,z^n) = \int
D_N(0,\xi,z)\prod_{j=1}^n\big(\|\xi^j\|^2 d\xi^j\big) d \xi$
\cite{BSZ1}.

The scaling asymptotics of the correlation functions then reduce to scaling asymptotics of the Bergman-\szego
kernel: In normal coordinates $\{z_j\}$ at $P_0\in M$ and in a `preferred' local frame for $L$, we have
\cite{BSZ1}:
\begin{eqnarray*}
  & & \frac{\pi^m}{N^m}\Pi_N(P_0+\frac{u}{\sqrtn},\frac{\theta}{N}; P_0+\frac{v}{\sqrtn},\frac{\phi}{N}) \\
  & \sim & e^{i(\theta-\phi)+u \cdot\bar{v} - \half(|u|^2 + |v|^2)} \textstyle \left[1+ b_{1}(u,v) N^{-\frac{1}{2}} +
  \cdots\right].
\end{eqnarray*}
To be precise, $\Pi_N$ is the natural lift of the kernel as an equivariant kernel on the boundary $\partial D^*$
of the unit (co-) disc bundle of $L^*$.  Note that $e^{i(\theta-\phi)+u \cdot\bar{v} - \half(|u|^2 + |v|^2)}$ is
the Bergman-\szego kernel of the Heisenberg group. These asymptotics use the Boutet de Monvel -Sjostrand
parametrix for the Bergman-\szego  kernel \cite{BS}, as applied in  \cite{Z4} to the  Fourier coefficients of the
kernel on powers of positive line bundles.

\subsection{Polynomials with fixed Newton polytope}
\vskip-5mm \hspace{5mm }

The well-known Bernstein-Kouchnirenko theorem states that the
number of simultaneous zeros of (a generic family of) $m$
polynomials with Newton polytope $P$ equals $m! Vol(P).$ Recall
that the Newton polytope $P_f$ of a polynomial is the convex hull
of its support $S_f = \{\alpha\in\Z^m: c_{\alpha} \not = 0\}$.
Using the homogenization map $f \to F$, the space of polynomials
$f$ whose  Newton polytope $P_f$ contained in $P$ may be
identified with  a subspace
\begin{equation} \label{SUBSPACE} H^0(\CP^m,
{\mathcal O}(p), P) = \{F \in  H^0(\CP^m, {\mathcal O}(p)): P_f
\subset P\}
\end{equation}
of $H^0(\CP^m, {\mathcal O}(p)).$

 The problem we address in this section is:
\begin{itemize}

\item {\bf Problem 3} How does the  Newton polytope
influence on the distribution of
 zeros of polynomials?

\end{itemize}

Again, one could ask the same question about $L^2$ mass, critical
points and so on and obtain a similar story. In \cite{SZ2} we
explore this influence in a statistical and asymptotic sense. The
main theme is that for each property of polynomials under study,
$P$ gives rise to  {\it classically allowed regions}  where the
behavior is the same as if no condition were placed on the
polynomials, and {\it classically forbidden regions} where the
behavior is exotic.

Let us define these terms.
 If $P \subset \R_+^m$ is a  convex integral polytope, then the  {\it classically allowed region\/} for polynomials in $H^0(\CP^m, {\mathcal
O}(p), P)$ is the set
$${\mathcal A}_P :=\mu_\Si\inv\left(\frac{1}{p}P^\circ\right) \subset \C^{*m}$$
(where $P^\circ$ denotes the interior of $P$), and the {\it
classically forbidden region\/} is its complement
$\C^{*m}\sm{\mathcal A}_P $. Here,
$\mu_\Si(z)=\left(\frac{|z_1|^2}{1+\|z\|^2},\dots,
\frac{|z_m|^2}{1+\|z\|^2} \right)\;$ is the moment map of $\CP^m.$

The result alluded to above is statistical. Since we view the
polytope $P$ of degree $p$ as placing a condition on the Gaussian
ensemble of $SU(P)$ polynomials of degree $p$,  we endow
$H^0(\CP^m, {\mathcal O}(p), P)$ with the   {\it conditional
probability measure\/} $\gamma_\de|_P$:
 \begin{equation} \label{CG} d \gamma_{\de|P} (s) =
\frac{1}{\pi^{\# P}}e^{-|\la|^2} d\la,\quad s = \sum_{\alpha \in
P} \la_{\alpha} \frac{z^{\alpha}}{\|z^{\alpha}\|}\;,
\end{equation} where the coefficients $\la_{\alpha}$ are again
independent complex Gaussian random variables with mean zero and
variance one.

Our simplest result concerns
  the  the expected density $\E_{|P} (Z_{f_1,
\dots, f_m})$  of the simultaneous zeros of $(f_1, \dots, f_m)$
chosen independently from $H^0(\CP^m, {\mathcal O}(p), P)$. It is
the measure on $\C^{*m}$ given by
\begin{eqnarray}\label{EZ0}
  \hspace*{-1cm} & & \E_{|P} (Z_{f_1, \dots, f_m})(U) \nonumber \\
  \hspace*{-1cm} & = & \int d\ga_{p|P}(f_1)\cdots\int d\ga_{p|P}(f_m)\;\big[\#\{z\in U:
f_1(z)=\cdots=f_m(z)=0\}\big]\;,
\end{eqnarray}
for  $U\subset \C^{*m}$, where the integrals are over $H^0(\CP^m,\ocal(p),P)$. We will determine the asymptotics
of the expected density as the polytope is dilated $P \to N P, N \in \N. $

\begin{theo}\label{probK} {\rm \cite{SZ2}}
Suppose that $P$ is a simple polytope in $\R^m$. Then, as $P$ is
dilated to  $NP$,
$$\frac{1}{(N\de)^m}\E_{|NP} (Z_{f_1, \dots, f_m})\to
\left\{\begin{array}{ll}\om_\FS^m
\ \ & \mbox{\rm on \ }{\mathcal A}_P \\[10pt]
0 & \mbox{\rm on \ }\C^{*m}\sm {\mathcal A}_P \end{array} \right.\
,$$ in the distribution sense; i.e., for any open
$U\subset\C^{*m}$, we have
$$\frac{1}{(N\de)^m}\E_{|NP}\big(\#\{z\in U:
f_1(z)=\cdots=f_m(z)=0\}\big)\to m!\vol_{\CP^m}(U\cap {\mathcal
A}_P )\;.$$
\end{theo}

There are also results for $k < m$ polynomials. The  distribution
of zeros is $\om_\FS^k$ in ${\mathcal A}_P$ as if there were no
constraint, while there is an exotic distribution in $\C^{*m}\sm
{\mathcal A}_P$ which depends on the exponentially decaying
asymptotics of
  the conditional Bergman- \szego
kernel
$$\Pi_{| NP} (z,w) = \sum_{\alpha \in NP} \frac{z^{\alpha}\overline{w}^{\alpha}}{||z^{\alpha}||_{FS} ||w^{\alpha}||_{FS} },$$
i.e. the orthogonal projection onto the subspace (\ref{SUBSPACE}).
It is obtained by sifting out terms in the (elementary) \szego
projector of $H^0(\CP^m, {\mathcal O}(p N))$ using the polytope
character $\chi_{NP}(e^{i \phi}) = \sum_{\alpha \in NP} e^{i
\langle \alpha, \phi \rangle}. $ To obtain asymptotics in the
forbidden region, we write  $\chi_{NP}(e^{i \phi}) = \int_{M_P}
\Pi_N^{M_P}(e^{i \phi}  w, w) dV(w)$, where $\Pi^{M_P}$ is
Bergman-\szego kernel of the toric variety $M_P$ associated to
$P$. We then make an explicit construction of $\Pi_N^{M_P}$ as a
complex oscillatory integral.  An alternative is to express
$\chi_{NP}$ as a Todd derivative of an exponential integral over
$P$ (following works of Khovanskii-Pukhlikov, Brion-Vergne and
Guillemin).  We thus obtain a complex oscillatory integral formula
for $\Pi_{|NP}(z,w)$. To obtain asymptotics in the forbidden
region we carefully deform the contour into the complex and apply
a complex stationary phase   method.

Although we only discuss expected behavior of zeros here, the
distribution of zeros is {\it self-averaging}: i.e., almost all
polynomials exhibit the expected behavior in an asymptotic sense.
We also expect similar results for critical points.

\section{Eigenfunctions} \label{section 1}\setzero

\vskip-5mm \hspace{5mm}

We now turn to the eigenvalue problem $\Delta_g \phi_{\nu}=
\lambda_{\nu}^2 \phi_{\nu},\; \langle \phi_{\nu}, \phi_{\nu'}
\rangle = \delta_{\nu \nu'}$ on a compact Riemannian manifold $(M,
g)$. We denote the $\lambda$-eigenspace by $V_{\lambda}.$ The role
of the \szego kernel is now played by the kernel $
E_\lambda(x,y)=\sum_{\lambda_\nu\le
\lambda}\phi_\nu(x)\overline{\phi_\nu(y)}$
 of the spectral projections.

\subsection{{\boldmath $L^p$} bounds}
\vskip-5mm \hspace{5mm }

Our first concern is with $L^p$ norms of $L^2$-normalized
eigenfunctions.
 We measure the growth rate of $L^p$ norms by
$ L^{p}(\lambda, g) = \sup_{\phi\in V_{\lambda}: ||\phi||_{L^2} =
1 }
 ||\phi||_{L^{p}}.  $
By the local Weyl law, $ E_{\lambda}(x,x)  = \sum_{\lambda_{\nu}
\leq \lambda} |\phi_{\nu}(x)|^2 = (2\pi)^{-n} \int_{|\xi|\le
\lambda} d\xi + O(\lambda^{n-1}),$ it follows that $
L^{\infty}(\lambda, g) = 0(\lambda^{\frac{n-1}{2}})$ on any
compact Riemannian manifold. This  bound, which is based entirely
on a local analysis, is sharp in the case of the standard round
sphere, $S^n$ or on any rotationally invariant metric on $S^2$,
but is far off in the case of flat tori. This motivates:

\begin{itemize}

\item {\bf Problem 5} For which $(M, g)$ is this estimate sharp? Which $(M, g)$
are extremal for growth rates of $||\phi_{\lambda}||_p$, both
maximal and minimal? What if $M$ has a boundary?  What is the
expected $L^p$ norm of a `random' $L^2$-normalized  polynomial or
eigenfunction?

\end{itemize}

In \cite{SZ}, we  give a  necessary condition for  {\it maximal}
eigenfunction growth: there  must exist a point $x \in M$ for
which the set $  {\mathcal L}_x = \{ \xi \in S^*_xM : \exists T:
\exp_x T \xi = x\} $ of directions of geodesic loops at $x$ has
positive surface measure.

\begin{theo}\label{eigcor} {\rm \cite{SZ}}
If  ${\mathcal L}_x$ has measure $0$ in $S^*_x M$ for every $x\in
M$ then
\begin{equation}\label{eig2}
L^p(\lambda, g)  = o( \lambda^{\delta(p)}), \, \, \, \,
p>\tfrac{2(n+1)}{n-1}\, \;\; \delta(p)=
\begin{cases}
n(\tfrac12-\tfrac1p)-\tfrac12, \quad \tfrac{2(n+1)}{n-1}\le p\le
\infty
\\
\tfrac{n-1}2(\tfrac12-\tfrac1p),\quad 2\le p\le
\tfrac{2(n+1)}{n-1}.
\end{cases}
\end{equation}
\end{theo}

The $L^p$-bounds  $O ( \lambda^{\delta(p)})$ were proved by Sogge
to hold for all $(M, g)$.

We further prove:

\begin{theo} \label{ANAL0} {\rm \cite{SZ}} (see also \cite{S}) Suppose that $(M,g)$ is:
\begin{itemize}

\item  Real analytic and that    $L^{\infty}(\lambda,g) = \Omega(\lambda^{(n-1)/2})$.
Then $(M,g)$ is a  $ Y^m_{\ell}$ -manifold, i.e. $\exists m$ such
that all geodesics    issuing from the point $m$ return to $m$ at
time  $\ell$. In particular, if dim $M$ = 2, then $M$ is
topologically a 2-sphere $S^2$ or a real projective plane $\R
P^2$.

\item  Generic. Then  $L^{\infty}(\lambda,g) = o (\lambda^{(n-1)/2})$.
\end{itemize}

 \end{theo}

\noindent Here,  $\Omega(\lambda^{\frac{n-1}2})$ means
$O(\lambda^{\frac{n-1}2})$ but not $o(\lambda^{\frac{n-1}2})$. The
generic result holds because
 ${\mathcal L}_x $ has measure $  0$ in $S^*_x M$ for all $ x \in M$
 for a residual set of metrics.

 In the case of random polynomials, or random combinations of
 eigenfunctions in short spectral intervals,  the almost sure growth of $L^{\infty}$
 norms is $O(\sqrt{\log N})$ while the $L^p$ norms for $p <
 \infty$ are bounded. This was proved by J. Vanderkam \cite{V} for $S^m$, Nonnenmacher-Voros  \cite{NV} for
 elliptic curves and Shiffman-Zelditch (to appear) for the general
 case using Levy concentration of measure estimates.

\subsection{Integrable case}
\vskip-5mm \hspace{5mm }

Results on  {\it minimal growth} have been obtained by J. A. Toth
and the author in the
 {\it quantum completely integrable} case,
where $\sqrt{\Delta} = P_1$ commutes with $n - 1$ first order
pseudodifferential operators $P_2, \dots, P_n \in \Psi^1(M)$ ($n$
= dim $M$)
 satisfying
$[P_i, P_j]=0$ and whose  symbols define a moment map ${\mathcal
P}:= (p_1, \dots, p_n)$ satisfying $dp_1\wedge dp_2 \wedge \dots
\wedge dp_n \not= 0$ on a dense open set $\Omega \subset T^*M -
0$. Since $\{p_i, p_j\} = 0$, the functions $p_1, \dots, p_n$
generate a homogeneous Hamiltonian $\R^n$-action whose  orbits
foliate $T^*M - 0$.   We refer to this foliation as the {\it
Liouville foliation}.

We consider the $L^p$ norms of the $L^2$-normalized joint
eigenfunctions $P_j \phi_{\lambda} = \lambda_j \phi_{\lambda}$.
The spectrum of $\Delta$ often has bounded multiplicity, so the
behaviour of joint eigenfunctions has implications for all
eigenfunctions.

\begin{theo} \label{LP} {\rm \cite{TZ1, TZ2}}   Suppose that the Laplacian $\Delta_g$ of   $(M, g)$ is quantum completely integrable
and that the joint eigenfunctions have uniformly bounded
$L^{\infty}$ norms. Then $(M, g)$ is a flat torus.
\end{theo}

This is  a kind of quantum analogue of the `Hopf conjecture'
(proved by Burago-Ivanov) that metrics on tori without conjugate
points are flat. In \cite{TZ2}, a quantitative improvement is
given under a further non-degeneracy assumption.  Unless $(M, g)$
is a flat torus, the Liouville foliation must possess a singular
leaf of dimension $ < n.$ Let $\ell$ denote the minimum dimension
of the leaves. We then construct a sequence of eigenfunctions
satisfying:
$$    \| \phi_{k} \|_{L^{\infty}} \geq C(\epsilon) \lambda_{k}^{\frac{n - \ell}{4} -
\epsilon},\;\;
 \| \phi_{k} \|_{L^{p}} \geq C(\epsilon) \lambda_{k}^{ \frac{(n - \ell) ( p - 2)}{ 4p } - \epsilon },
 \,\,\,(  2 < p)$$
for any $\epsilon > 0$. It is easy to construct examples were
$\ell = n - 1$, but it seems plausible that in `many' cases $\ell
= 1$. To investigate this, one would study the boundary faces of
the image ${\mathcal P}(T^*M - 0)$ of $T*M - 0$ under a
homogeneous moment map. For a related study in the case of torus
actions, see Lerman-Shirokova \cite{LS}.

\subsection{Quantum ergodicity}
\vskip-5mm \hspace{5mm }

Quantum ergodicity is concerned with the sums ($A \in \Psi^0(M)$):
\begin{equation} \label{QE} S_p(\lambda) = \sum_{\nu: \lambda_{\nu} \leq \lambda}
 |\langle A \phi_{\nu}, \phi_{\nu} \rangle - \omega(A)|^p,\;\;
 \omega(A) = \frac{1}{Vol(S^*M)} \int_{S^*M} \sigma_A d \mu. \end{equation}

In work of A.I. Schnirelman \cite{LSchn}, Colin de Verdiere and
the author \cite{Z2}, it is shown that $S_p(\lambda) =
o(N(\lambda))$ if $G^t$ is ergodic. In the author's view
\cite{Z2}, this is best viewed as a convexity theorem.  We mention
briefly some new results.

In work  of Gerard-Leichtnam \cite{GL} and Zelditch-Zworski
\cite{ZZw}, the ergodicity result was  extended  to domains
$\partial \Omega$ with piecewise smooth boundary and ergodic
billiard flow. Since the billiard map on $B^* \partial \Omega$ is
ergodic whenever the billiard flow is, suitable boundary values of
ergodic eigenfunctions (e.g. $\phi_k|_{\partial \Omega}$ in the
Neumann case or  $\partial_{\nu} \phi_k|_{\partial \Omega}$ in the
Dirichlet case) should also have the ergodic property. This was
conjectured by S. Ozawa in 1993. A proof is given in our  work
with A. Hassell \cite{HZ} for convex piecewise smooth domains with
ergodic billiards (in the case of domains with Lipschitz normal
and with Dirichlet boundary conditions, this had earlier been
proved in \cite{GL} by a different method).

Little is known about the rate of decay. For sufficiently chaotic
systems (satisfying the central limit theorem) one can get the
tiny improvement  $S_p(\lambda) = O(N(\lambda)/(\log \lambda)^p))$
\cite{Z3}.  The asymptotics $S_2(\lambda) \sim B(A) \lambda$ have
recently been obtained by Luo-Sarnak \cite{LuS} for Hecke
eigenfunctions of the modular group, exploiting the connections
with $L$-functions. These asymptotics (though not the coefficient)
are predicted by the random polynomial model.  Other strong bounds
 in the arithmetic case were obtained by Kurlberg-Rudnick for eigenfunctions of certain
quantized torus automorphisms \cite{KZ}. Bourgain-Lindenstrauss
\cite{BL}   and Wolpert \cite{W} have developed the 'non-scarring'
result  of \cite{RS} to give entropy estimates of possible quantum
limit measures in arithmetic cases.

A natural problem is the converse:
\begin{itemize}

\item {\bf Problem 6} What can be said of the dynamics if $S_p(\lambda) = o(N(\lambda))$?
Does quantum ergodicity imply classical ergodicity?
\end{itemize}

It is known that classical ergodicity is equivalent to this bound
plus estimates on off-diagonal terms \cite{Su}. The existence of
KAM quasimodes (due to  Lazutkin \cite{LSchn}, Colin de Verdiere
\cite{CV}, and Popov \cite{P}) makes it very likely that KAM
systems are {\it not} quantum ergodic, nor are $(M, g)$ which have
stable elliptic orbits.

A further  problem which may be accessible is:
\begin{itemize}

\item {\bf Problem 7}  How are the nodal sets $\{\phi_{\nu} = 0\}$ distributed in
the limit $\nu \to \infty$?
\end{itemize}

In \cite{NV} (for elliptic curves) and \cite{SZ2} (general \kahler
manifolds) it is proved that the {\it complex zeros} of quantum
ergodic eigenfunctions become uniformly distributed relative to
the volume form. Can one prove an analogue for the real zeros?

\label{lastpage}

\end{document}